# An Optimization Framework for the Time-Dependent Electric Vehicle Routing Problem with Shared Mobility: A Step Toward Smart Cities

**Alireza Yazdiani**
Department of Civil and Environmental Engineering
Cornell University, Ithaca, USA, 14853
Email: ay373@cornell.edu
https://orcid.org/0000-0002-4135-5209

**Shayan Bafandkar**
The Grainger College of Engineering, Department of Civil and Environmental Engineering, University of Illinois at Urbana-Champaign, Urbana, IL 61801, USA
Email: sbafan2@illinois.edu
https://orcid.org/0009-0009-8172-5751

**Amir Elmi**
School of Civil Engineering, The University of Sydney, Sydney, Australia
Email: amir.elmi@sydney.edu.au
https://orcid.org/0000-0003-0127-7573

**Yousef Shafahi**
Department of Civil Engineering
Sharif University of Technology, Tehran, Iran, 1458889694
Email: shafahi@sharif.edu
https://orcid.org/0000-0003-4267-4348

**Corresponding author**
Alireza Yazdiani, ay373@cornell.edu

**Abstract**

This paper aims to introduce a mathematical model to solve the time-dependent electric vehicles routing problem in shared travels. Shared mobility has gained significance recently due to its contribution to the alleviation of traffic congestion and air pollution. In this study, a MIP model has been developed and solved for a toy problem using the CPLEX solver to indicate the efficiency of shared mobility compared to private mode. We have considered practical constraints, such as considering the traffic congestion throughout the day by assigning a time-dependent step function to the network's links altering travel times in the peak and off-peak hours, nonlinear charging function, vehicles' queue at charging stations, partial charging possibility, different charging infrastructures, and passengers' desired timewindows, to make the results of the model more applicable in the real world. Among the important results of this presented model according to the solved example, we can mention the reduction of mileage per vehicle in personal mode from 56.42 km to 46.83 km in shared mode while the total travel time per vehicle has only increased slightly from $84.63\frac{min}{veh}$ to $89.32\frac{min}{veh}$, respectively. This suggests that the typical problems related to private travel can be offset by shared travel without losing the comfort and convenience of the former.

**Keywords**

### 1. Introduction

In today's world, a high level of urbanism has led to an increase in travel demand in urban areas, most of which (40% to 83%) is made using privately owned vehicles [1]. Although the private mode of transportation brings about more convenience and comfort for passengers, it places problematic consequences on society due to their low utility, need for parking space, and environmental issues. In fact, Litman et al. [2] indicated in their study of private cars' utility in Austin, US, that they are out of service 23 hours a day on average. To address this issue, several programs have been adopted by urban decision-makers, one of which is shifting a proportion of the daily travels to more sustainable modes of transportation, including shared mobility, which is in the scope of this study.

Shared mobility, as a sustainable mode of transportation, has made a huge impact on private car ownership, average transportation cost, greenhouse gas emission, and traffic congestion [3]. Shared mobility, including different services such as car sharing, ride sharing, and bike sharing, has gained significance recently and will soon replace traditional modes of public transportation [4]. During the last decade, shared travels made up 50% of private travel demand in the San Francisco Bay [5]. Furthermore, Alonso-Mora et al. [6] showed that shared mobility has reduced taxi rides by 75%, and Talebpour and Mahmassani [7] suggested that traffic congestion could be heavily reduced because of shared mobility. A portion of the vehicle fleet in this mode consists of electric vehicles. These vehicles enhance the benefits of shared mobility by not only improving travel safety but also contributing positively to the environment [8]. They can be used in a dial-a-ride transit system which can provide regular shared travels' flexibility while maintaining taxi rides' features such as low accessibility cost, high capacity, and booking possibility. However, use of electric vehicles in dial-a-ride transit is still facing ongoing challenges due to the limited fleet size, their charging process that need to be addressed by redesigning the transit system enabling it to serve high demand on a large scale.

We propose an optimization model for the time-dependent electric vehicles routing problem in the shared travel mode to determine the trip plan of each electric vehicle in the network according to a pre-determined static demand. The trip demand is considered static and passengers have determined their origin and destination along with their desired pick-up from origin and arrival at destination time window. Our research contributes to the dial-a-ride transit system literature in the following ways: (i) we introduce a time-dependent formulation capturing network's traffic congestion using a step speed-time function; (ii) we consider different waiting times at the charging stations throughout the day and several charging infrastructures with different charging powers enabling vehicles to manage their battery level according to the demand, waiting time, and their battery levels; and (iii) we provide a shared-trip plan that could reduce the operational cost of relocating electric vehicles and traffic congestion in urban areas. The following is divided into four sections. In section 2, previous related studies are briefly reviewed. Then, the optimization model is introduced and solved using the CPLEX solver for a toy problem to indicate its efficiency in sections 3 and 4, respectively. Finally, conclusions are made in section 5.

## 2. Literature Review

Passenger trip business models, such as those used in shared mobility services and public transportation, are often formulated as optimization models aimed at enhancing utility features, including minimizing passenger waiting times and improving user interfaces between drivers and passengers in carpooling, vanpooling, and peer-to-peer (P2P) ride-sharing scenarios ([9], [10]). One of the main focuses of the mentioned models is the vehicle routing problem (VRP). Toth and Vigo [11] provided a general definition of the classic VRP. The constraints of classic VRP focus on routing the vehicle fleet to meet passenger and freight demand according to a defined objective function. Almeida et al. [12] presented a more detailed explanation of VRP by dividing it into several categories based on the available operations research literature. Vehicle Routing Problem (VRP) types include Capacitated VRP (fleet capacity), VRP with time windows (service within a time window), VRP with multiple trips (multiple trips before depot return), and Open VRP (no depot return constraint). Dial-a-ride transport (DART) or demand responsive transport (DRT) includes a pick-up time window and trip time limit. Time-dependent VRP with time windows (TDVRPTW) considers time-varying trip times, and Dynamic Vehicle Routing uses real-time data for rerouting. Our study focuses on Time-dependent Dial-a-ride Problem (TD-DARP), a blend of Capacitated VRP, DRT, and TDVRPTW.

DARP is a complex optimization problem that finds its utility in on-demand transportation systems. Cordeau and Laporte [13] attribute the complexity of DARP to various cost weight factors and considerations related to passengers' utility. The complexity is further compounded when considering the time dependency of trip times (TD-DARP). This necessitates the minimization of passengers' trip length according to the current time in the model, along with passengers' utility. Several strategies are available in the literature to simplify the TD-DARP. For instance, Malandraki and Daskin [14] assigned a step-time function to links' travel times and solved the TD-DARP for a pre-determined passenger demand using the nearest-neighbor metaheuristic method. Moreover, the electrification of the vehicle fleet adds another layer of complexity to TD-DARP as these vehicles have relatively high charging times and the charging function should be involved in the formulation. The existing literature on resolving the electric vehicle routing problem is somewhat limited. Zhang et al. [15] formulated an optimization model to minimize the energy consumption of the electric vehicle fleet in the routing process due to the limited battery capacity, passenger's range anxiety, and environmental concerns. They indicated that the ant colony heuristic method is suitable for solving the problem. Keskin et al. [16] worked on strategies for partial recharging within soft time windows when electric vehicles have to wait in a queue each time they want to recharge. Five charging stations' queues are defined to consider the altering hourly traffic flow. Their results indicate that although the total traveled distance does not change significantly, the queue time has an important effect on the vehicles' routes. Montoya et al. [17] solved the electric VRP with a nonlinear charging function by dividing the function into several sections and making linear approximations of the function at each section. Hiermann et al. [18] tackled the VRP with a heterogeneous fleet composed of both electric and regular vehicles by developing an adaptive large neighborhood search method.

According to the authors' knowledge and based on the conducted literature review on TD-DARP, none of the previous studies in this area considered traffic congestion through out the network, several charging stations with different charging powers in the network, the possibility of partial charging, and an electric vehicle fleet with a nonlinear charging function simultaneously. Therefore, the novelty of this paper lies in its focus on introducing a novel methodology to solve the TD-DARP with static definite demand and an electric vehicle fleet while considering 50 constraints in the formulation, making the results more applicable in the real world. In the end, it is implemented on a toy problem as a case study, and the derived results are thoroughly explained. The details are provided in the following.

### 3. Methodology

Our model aims to minimize a weighted travel cost function for the electric vehicle TD-DARP while meeting static demand. Passenger demand includes passenger count, pick-up/drop-off coordinates, and desired time windows. The network is represented as a graph $G = (V, A)$, with nodes (V) as pick-up/drop-off points, charging stations, and depot. The set (A) is a Euclidean distance matrix between nodes. Each demand is represented as a (i, n+i) pair, where 'i' is the pick-up point, 'n' is total demands, and 'n+i' is the drop-off point.

We centralize all trip demands at each traffic zone's core. Multiple passengers can board at the pick-up point if they share the same destination, subject to vehicle capacity. Electric vehicles must maintain a minimum battery level for battery health. If the battery level falls below the threshold, they must find the most cost-effective charging station considering travel, waiting, and charging times. A vehicle may prefer a distant station with a short queue over a nearby station with a long queue.

Vehicle energy consumption depends on the traveled distance, allowing for an energy consumption matrix calculation using the Euclidean distance matrix. To account for traffic congestion, we employed the method of Malandraki and Daskin [14], considering varying travel speeds for vehicles at different times of the day. A step-time function is assigned to each link based on the distance and speed, with a sample shown in Figure 1.

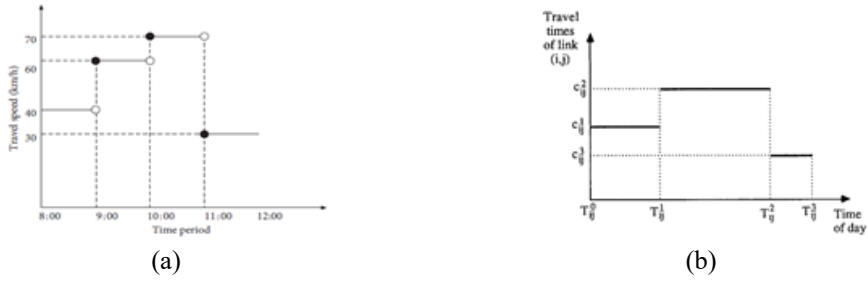

(a)                                                    (b)

**Figure 1-** The step-speed (a) and step-time (b) functions introduced by Malandraki and Daskin [14]

We adopted Montoya et al. [17]'s nonlinear charging function and three different charging infrastructures. Therefore, the charging rate at each station varies based on the vehicle's battery level, linear approximation of the charging function, and the station's charging power. These two concepts are depicted in Figures 2 and 3, respectively.

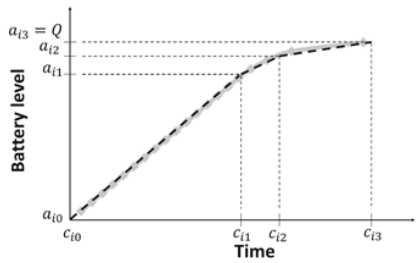

**Figure 2-** The nonlinear charging fnction with linear approximations introduced by Montoya et al. [17]

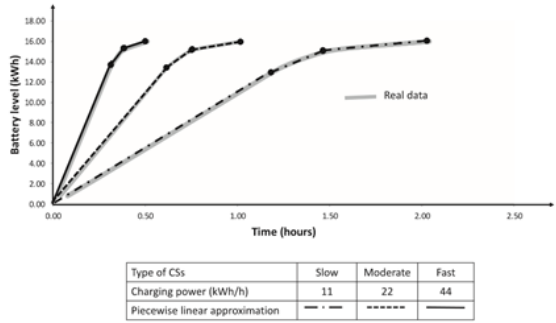

**Figure 3-** Three charging stations with different powers in the tudy of Montoya et al. [17]

**Table1**

The breakpoints' coordinates at the linear approximations of the charging function at different stations introduced in the Montoya et al. [17]'s study:

| Variable | Slow charging station (11kWh/h) | Moderate charging station (22kWh/h) | Fast charging station (44kWh/h) |
|---|---|---|---|
| $a_{i,r}^k$: The battery level of vehicle k upon reaching | $a_{1,0} = 0$ $a_{1,1} = 12.2$ | $a_{2,0} = 0$ $a_{2,1} = 13.1$ | $a_{3,0} = 0$ $a_{3,1} = 13.9$ |

| break point $r \in B$ while charging at station $i \in F$ | $a_{1,2} = 14.9$ $a_{1,3} = 16$ | $a_{2,2} = 15.2$ $a_{2,3} = 16$ | $a_{3,2} = 15.3$ $a_{3,3} = 16$ |
|---|---|---|---|
| $c_{i,r}^k$: The charging time (min) of vehicle $k \in K$ to reach break point $r \in B$ while charging at station $i \in F$ | $c_{1,0} = 0$ $c_{1,1} = 72$ $c_{1,2} = 85$ $c_{1,3} = 125$ | $c_{2,0} = 0$ $c_{2,1} = 32$ $c_{2,2} = 41$ $c_{2,3} = 58$ | $c_{3,0} = 0$ $c_{3,1} = 15$ $c_{3,2} = 21$ $c_{3,3} = 28$ |

Queue scenarios at each station are designed based on the study by Keskin et al. [16]. According to this methodology, the depot's charging station is considered to be without a queue, while other charging stations have queues that vary with time. All the scenarios under consideration are detailed in Table 2. Eventually, it has been assumed that in a dial-a-ride shared trip system, passengers consent to travel with other passengers in the vehicle.

**Table2**

Keskin et al. [16]'s queue waiting times at charging stations (min):

| Scenario | Queue length | Morning | Noon | Late afternoon | Evening |
|---|---|---|---|---|---|
| No waiting | - | 0 | 0 | 0 | 0 |
| Time-independent | Short | 10 | 10 | 10 | 10 |
| | Long | 20 | 20 | 20 | 20 |
| Time-dependent/Smooth transitions | Short | 10 | 20 | 10 | 5 |
| | Long | 20 | 40 | 20 | 10 |
| Time-dependent/Steep transitions | Short | 20 | 50 | 25 | 5 |
| | Long | 40 | 100 | 50 | 10 |

A 4-index formulation for the electric TD-DARP, building upon the DARP formulation provided by Cordeau [19] and Bongiovanni [20], is developed allowing the model to account for time-dependency in the routing process. A binary decision variable $x_{i,j}^{k,m}$ indicates whether vehicle k sequentially visits nodes i and j in the $m^{th}$ time period. If vehicle k ($k \in K$) departs from node $i \in v$ heading to node $j \in v$ in the $m^{th}$ time period ($m \in M$), $x_{i,j}^{k,m} = 1$; otherwise, $x_{i,j}^{k,m} = 0$. The decision variable $U_i^k$ denotes the time when vehicle k begins its service at node i. The model's sets, parameters, and auxiliary variables are mentioned in Table 3, respectively.

**Table3**

Variables and parameters of the optimization model:

| Notation | Definition |
|---|---|
| P = $\{1, \dots, n\}$ | Set of pick-up points |
| D = $\{n + 1, \dots, 2n\}$ | Set of drop-off points |
| K = $\{1, \dots, k\}$ | Set of vehicle fleet |
| M = $\{1, \dots, m\}$ | Set of time periods |
| $F = \{1, \dots, f\}$ | Set of charging stations |
| $V = 0 \cup P \cup D \cup F \cup 2n+1$ | Set of all possible stop locations |
| $B = \{1, \dots, b\}$ | Set of break points at the linear approximations of the charging function |
| $0$ | Initial depot node |
| $2n + 1$ | Final depot node (may be different from the initial depot) |
| $Q^k$ | Vehicle k's capacity parameter |
| $T^k$ | The allowable amount of time that vehicle k could be used |
| $early_i$ | Earliest time that the service could begin at node $i \in v$ |
| $late_i$ | Latest time that the service could begin at node $i \in v$ |
| $q_i$ | Variation in the number of embarked passengers at node $i \in v$ |
| $dstop_i$ | Service (stop) time at node $i \in v$ |
| $B_{cap}^k$ | The battery capacity of vehicle $k \in K$ |
| $b^{m-1}$ | The beginning of time period $m \in M$ |
| $V_{i,j}^m$ | The vehicle's speed on link ij ($\forall i \in V, j \in V, i \neq j$) in time period $m \in M$ (represented on the vertical axis of Figure 1 (a)) |
| $t_{i,j}^m$ | The travel time of link ij ($\forall i \in V, j \in V, i \neq j$) in time period $m \in M$ |
| $queue_i^m$ | Waiting time at charging station $i \in F$ in time period $m \in M$ |
| $E0^k$ | The battery level of vehicle k when leaving initial depot node |
| $B_{Min}$ | Minimum allowable battery level indicated by the proportion of the battery capacity |
| $B_{Max}$ | Maximum allowable battery level indicated by the proportion of the battery capacity |
| $C_{Min}$ | The minimum amount of charging when entering a charging station indicated by the proportion of the battery capacity |
| $e_{ij}$ | Energy consumed while traveling from node $i \in v$ to node $f \in v$ (distance dependent) |
| $a_{i,r}^k$ | The battery level of vehicle k upon reaching break point $r \in B$ while charging at station $i \in F$ (computed according to Table 1) |
| $c_{i,r}^k$ | The charging time of vehicle $k \in K$ to reach break point $r \in B$ while charging at station $i \in F$ (computed according to Table 1) |
| $M_t$ | A very large number used in the constraints related to trip times |
| $M_p$ | A very large number used in the constraints related to passengers |
| $M_B$ | A very large number used in the constraints related to vehicles' batteries |
| $E1_i^k$ | Auxiliary variable indicating the battery level of vehicle $k \in K$ when reaching node $i \in v$ |
| $E2_i^k$ | Auxiliary variable indicating the battery level of vehicle $k \in K$ when departing from node $i \in v$ |
| $w_i^k$ | The number of passengers in vehicle $k \in K$ at the time of its departure from node $i \in v$ |
| $v_i^k$ | The amount of time vehicle $k \in K$ spends at charging station $i \in F$ |
| $twv_i^1$ | The amount of earliness of arrival at location $i \in P \cup D$ compared to the lower limit of the time window announced by passengers at these points |
| $twv_i^2$ | The amount of lateness of arrival at location $i \in P \cup D$ compared to the upper limit of the time window announced by passengers at these points |
| $twv_k^3$ | The amount of lateness of arrival at the final depot node $i \in P \cup D$ compared to the allowable amount of time that vehicle $k \in K$ could be used ($T^k$) |
| $g1_i^k$ | The charging time assigned to vehicle $k \in K$ upon arrival at charging station $i \in F$ according to the linear approximation of the charging function |
| $g2_i^k$ | The charging time assigned to vehicle $k \in K$ when leaving charging station $i \in F$ according to the linear approximation of the charging function |
| $h1_{i,r}^k$ | Break point $r(\in B)$'s factor in the linear approximation of the charging function when vehicle $k \in K$ arrives at charging station $i \in F$ $\left(0 \leq h1_{i,r}^k \leq 1\right)$ |
| $h2_{i,r}^k$ | Break point $r(\in B)$'s factor in the linear approximation of the charging function when vehicle $k \in K$ leaves charging station $i \in F$ $\left(0 \leq hw_{i,r}^k \leq 1\right)$ |
| $y_{i,r}^k$ | A binary auxiliary variable indicating the linear approximation of the charging function upon arrival at charging station. If the battery level of vehicle k ($k \in K$) is found between $a_{i,r}^k$ and $a_{i,r-1}^k$ $(r \in B \backslash \{0\})$ when entering charging station $i \in F$, $y_{i,r}^k = 1$; otherwise, $y_{i,r}^k = 0$. |
| $z_{i,r}^k$ | A binary auxiliary variable indicating the linear approximation of the charging function when leaving charging station. If the battery level of vehicle k ($k \in K$) is found between $a_{i,r}^k$ and $a_{i,r-1}^k$ $(r \in B \backslash \{0\})$ when leaving charging station $i \in F$, $z_{i,r}^k = 1$; otherwise, $z_{i,r}^k = 0$. |

The objective function and all of its constraints are mentioned in Table 4. The description of each equation is provided.

**Table 4**

Objective function and constraints of the routing optimization model:

| Equation | No. |
|---|---|
| $Min \ Z = w_1 \sum_{k \in k} \sum_{i \in v} \sum_{j \in v, (j \neq i)} \sum_{m \in M} t_{i,j}^m \times x_{i,j}^{k,m} + w_2 \sum_{i \in P \cup D} q_i (twv_i^1 + twv_i^2) + w_3 \sum_{k \in K} twv_k^3$ | 1 |

$$\sum_{k \in K} \sum_{j \in P \cup D, (j \neq i)} \sum_{m \in M} x_{i,j}^{k,m} = 1; \ \forall i \in P \tag{2}$$

$$\sum_{j \in P \cup F \cup 2n+1} \sum_{m \in M} x_{0,j}^{k,m} = 1; \ \forall k \in K \tag{3}$$

$$\sum_{i \in D \cup F \cup 0} \sum_{m \in M} x_{i,2n+1}^{k,m} = 1; \ \forall k \in K \tag{4}$$

$$\sum_{j \in P \cup D, j \neq i} \sum_{m \in M} x_{i,j}^{k,m} - \sum_{j \in P \cup D, j \neq i} \sum_{m \in M} x_{j,n+i}^{k,m} = 0; \ \forall i \in P, k \in K \tag{5}$$

$$\sum_{j \in V, j \neq i} \sum_{m \in M} x_{i,j}^{k,m} - \sum_{j \in V, j \neq i} \sum_{m \in M} x_{j,i}^{k,m} = 0; \ \forall i \in P \cup D \cup F, k \in K \tag{6}$$

$$U_i^k + t_{stop}^i + \sum_{m \in M} t_{i,j}^m \times x_{i,j}^{k,m} - M_t \times \left(1 - \sum_{m \in M} x_{i,j}^{k,m}\right) \leq U_j^k; \ \forall i \in V, j \in V, k \in K, i \neq j \tag{7}$$

$$early_i - U_i^k \leq twv_i^1; \ 0 \leq twv_i^1; \ \forall i \in P \cup D, k \in K \tag{8}$$

$$U_i^k - late_i \leq twv_i^2; 0 \leq twv_i^2; \ \forall i \in P \cup D, k \in K \tag{9}$$

$$U_i^k + dstop_i + \sum_{m \in M} t_{i,n+i}^m \times x_{i,n+i}^{k,m} \leq U_{n+i}^k; \ \forall i \in P, k \in K \tag{10}$$

$$w_i^k + q_j - M_p \left(1 - \sum_{m \in M} x_{i,j}^{k,m}\right) \leq w_j^k; \ \forall i \in V, j \in V, k \in K, i \neq j \tag{11}$$

$$w_i^k + q_j + M_p \left(1 - \sum_{m \in M} x_{i,j}^{k,m}\right) \geq w_j^k; \ \forall i \in V, j \in V, k \in K, i \neq j \tag{12}$$

$$w_i^k = 0; \ \forall i \in 0 \cup F \cup 2n+1, k \in K \tag{13}$$

$$w_i^k \geq max(0, q_i); \ \forall i \in P \cup D, k \in K \tag{14}$$

$$w_i^k \leq min(Q^k, Q^k + q_i); \ \forall i \in P \cup D, k \in K \tag{15}$$

$$E1_i^k = E0^k; \ \forall k \in K \tag{16}$$

$$E1_i^k \geq B_{Min} * B_{cap}^k * \sum_{j \in V, j \neq i} \sum_{m \in M} x_{i,j}^{k,m}; \ \forall k \in K, i \in V \tag{17}$$

$$E2_i^k \leq B_{Max} * B_{cap}^k * \sum_{j \in V, j \neq i} \sum_{m \in M} x_{i,j}^{k,m}; \ \forall k \in K, i \in F \tag{18}$$

$$E2_i^k \geq C_{Min} * B_{cap}^k * \sum_{j \in V, j \neq i} \sum_{m \in M} x_{j,i}^{k,m}; \ \forall k \in K, i \in F \tag{19}$$

$$E1_i^k \leq E2_i^k; \ \forall k \in K, i \in F \tag{20}$$

$$E1_i^k - e_{i,j} - M_B * \left(1 - \sum_{m \in M} x_{i,j}^{k,m}\right) \leq E1_j^k; \ \forall k \in K, i \in V \backslash F, j \in V \backslash 0, i \neq j \tag{21}$$

$$E1_i^k - e_{i,j} + M_B * \left(1 - \sum_{m \in M} x_{i,j}^{k,m}\right) \geq E1_j^k; \ \forall k \in K, i \in V \backslash F, j \in V \backslash 0, i \neq j \tag{22}$$

$$E2_i^k - e_{i,j} - M_B * \left(1 - \sum_{m \in M} x_{i,j}^{k,m}\right) \leq E1_j^k; \ \forall k \in K, i \in F, j \in P \cup F \cup 2n+1, i \neq j \tag{23}$$

$$E2_i^k - e_{i,j} + M_B * \left(1 - \sum_{m \in M} x_{i,j}^{k,m}\right) \geq E1_j^k; \ \forall k \in K, i \in F, j \in P \cup F \cup 2n+1, i \neq j \tag{24}$$

$$E1_i^k = \sum_{r \in B} h1_{i,r}^k a_{i,r}^k; \ \forall k \in K, i \in F \tag{25}$$

$$g1_i^k = \sum_{r \in B} h1_{i,r}^k c_{i,r}^k; \ \forall k \in K, i \in F \tag{26}$$

$$\sum_{r \in B \backslash \{0\}} y_{i,r}^k = \sum_{j \in V} \sum_{m \in M} x_{i,j}^{k,m}; \ \forall k \in K, i \in F \tag{27}$$

$$\sum_{r \in B} h1_{i,r}^k = \sum_{r \in B \setminus \{0\}} y_{i,r}^k ; \ \forall k \in K, i \in F \tag{28}$$

$$h1_{i,0}^k \leq y_{i,1}^k; \ \forall k \in K, i \in F \tag{29}$$

$$h1_{i,r}^k \leq y_{i,r}^k + y_{i,r+1}^k; \ \forall k \in K, i \in F, r \in B \setminus \{0, b\} \tag{30}$$

$$h1_{i,b}^k \leq y_{i,b}^k; \ k \in K, i \in F \tag{31}$$

$$E2_i^k = \sum_{r \in B} h2_{i,r}^k \, a_{i,r}^k; \ \forall k \in K, i \in F \tag{32}$$

$$g2_i^k = \sum_{r \in B} h2_{i,r}^k \, c_{i,r}^k; \ \forall k \in K, i \in F \tag{33}$$

$$\sum_{r \in B \setminus \{0\}} z_{i,r}^k = \sum_{j \in V \setminus 2n+1} \sum_{m \in M} x_{i,j}^{k,m}; \ \forall k \in K, i \in F \tag{34}$$

$$\sum_{r \in B} h2_{i,r}^k = \sum_{r \in B \setminus \{0\}} z_{i,r}^k ; \ \forall k \in K, i \in F \tag{35}$$

$$h2_{i,0}^k \leq z_{i,1}^k; \ \forall k \in K, i \in F \tag{36}$$

$$h2_{i,r}^k \leq z_{i,r}^k + z_{i,r+1}^k; \ \forall k \in K, i \in F, r \in B \setminus \{0, b\} \tag{37}$$

$$h2_{i,b}^k \leq z_{i,b}^k; \ \forall k \in K, i \in F \tag{38}$$

$$v_i^k = g2_i^k - g1_i^k; \ \forall k \in K, i \in F \tag{39}$$

$$v_i^k \leq U_j^k - \left( U_i^k + \sum_{m \in M} t_{i,j}^m * x_{i,j}^{k,m} + dstop_i + queue_i^m - M_t \left( 1 - \sum_{m \in M} x_{i,j}^{k,m} \right) \right); \ \forall i \in F, j \in V, k \in K, \ i \neq j \tag{40}$$

$$\sum_{k \in K} \sum_{i \in V \setminus 2n+1} \sum_{m \in M} x_{i,j}^{k,m} \leq 1; \ \forall j \in F \tag{41}$$

$$U_{2n+1}^k - T^k \leq twv_k^3; 0 \leq twv_k^3; \ \forall k \in K \tag{42}$$

$$U_i^k + dstop(i) \geq b_{i,j}^m * x_{i,j}^{k,m}; \ \forall i, j \in V, k \in K, m \in M \tag{43}$$

$$U_i^k + dstop(i) \leq b_{i,j}^{m+1} + M_t \times \left( 1 - x_{i,j}^{k,m} \right); \ \forall i, j \in V, k \in K, m \in M \tag{44}$$

$$h1_{i,r}^k \geq 0, \ h2_{i,r}^k \geq 0; \ \forall k \in K, i \in F, r \in B \tag{45}$$

$$h1_{i,r}^k \leq 1, \ h2_{i,r}^k \leq 1; \ \forall k \in K, i \in F, r \in B \tag{46}$$

$$g1_i^k \geq 0, \ g2_i^k \geq 0; \ \forall k \in K, i \in F \tag{47}$$

$$x_{i,j}^k \in \{0,1\}; \ \forall k \in K, i \in V, j \in V \tag{48}$$

$$E1_i^k \geq 0; \ \forall k \in K, i \in V \tag{49}$$

$$E2_i^k \geq 0; \ \forall k \in K, i \in F \tag{50}$$

$$v_i^k \geq 0; \ \forall k \in K, i \in F \tag{51}$$

$$y_{i,r}^k \in \{0,1\}, \ z_{i,r}^k \in \{0,1\}; \ \forall k \in K, i \in F, r \in B \setminus \{0\} \tag{52}$$

Equation (1) depicts the model's objective function, which consists of three terms. These terms represent the total travel time of vehicles considering traffic congestion, passengers' disutility caused by deviation from their desired timetables, and a delay penalty for late arrivals at the final depot, respectively. It's worth noting that if we were dealing with a heterogeneous fleet comprising both electric and regular vehicles, a fourth term could be added to account for the air pollution caused by regular vehicles.

Equation (2) ensures that each trip demand is served only once and by one vehicle throughout the day. Equations (3) and (4) denote that each vehicle commences its trip from the initial depot and ends it at the final depot. The number of passengers present at each station to embark on the train is determined by Equations (5) and (6).

Passenger demand is addressed in Equation (5), where vehicle k is compelled to visit node 'n+i' in subsequent time periods if it meets the passenger demand at node 'i' in time period 'm'. Equation (6) captures the flow balance at each node. A lower bound for the commencement of service at node

j∈V, which is visited immediately after node i∈V, is determined by Equation (7). The penalties for earliness and delay relative to passengers' desired time windows are computed in Equations (8) and (9). Equation (10) ensures that each passenger's drop-off point is visited after their pick-up point. The fluctuation in the number of passengers on board each vehicle at each node is depicted by Equations (11) and (12). Equation (13) captures another aspect of passenger utility, ensuring that vehicles are not occupied by passengers at charging stations and depots. The capacity of the vehicle is addressed through Equations (14) and (15).

Battery health preservation is considered in Equations (17) and (18), which determine the lower and upper allowable limits of the battery level. The minimum amount of charging at the stations is depicted by Equation (19). Equation (20) represents a lower bound for the battery level of vehicles leaving charging stations. Equations (21) and (22) encapsulate the energy consumed during travel from each origin node (except charging stations) to each destination node (except the initial depot). Equations (23) and (24) show the battery level of vehicles traveling from charging stations to their origins (except the initial depot).

Equations (25) to (31) determine the battery level and its corresponding charging time for a vehicle entering a charging station, according to the linear approximation of the charging function. Subsequently, Equations (32) to (38) determine the battery level and its corresponding charging time for a vehicle leaving a charging station, according to the linear approximation of the charging function. The amount of time each vehicle spends at each charging station and an upper limit for this variable are calculated using Equations (39) and (40), respectively. Equation (41) limits the number of visits to each charging station to one. The third term (final depot lateness penalty) is calculated by Equation (42)[1]. 

Equations (43) and (44) indicate the time period within which a link is traversed by a vehicle. Equations (45) to (52) are non-negativity and variable restriction constraints, setting feasibility intervals for the variables.

The introduced formulation is an integer programming (IP) model that could be solved by commercial optimization software such as GAMS. This study makes use of the CPLEX solver that relies on several algorithms including enumeration techniques such as the branch and bound approach as well as cutting-plane techniques to find the optimum solution, and has been widely used in other studies [21]. In the next section, the optimization model and algorithm introduced has been implemented on a toy problem to validate its competency.

## 4. Case Study

We designed a small-scale network for which the formulation introduced in Section 3 could be implemented and solved using the CPLEX solver. The network is a $20 \times 20 \ km^2$ square area with a daily demand of 5. The day is divided into 5 time periods, and the vehicles' travel speed changes across all links according to the step function depicted in Figure 4.

---

[1] As we are trying to minimize the objective function, this equation is bounded if the non-negativity of $twv_k^3$ is satisfied.

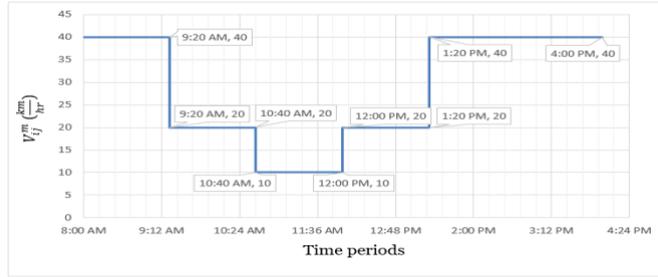

**Figure 4-** The step-speed function adopted in our study and introduced by Malandraki and Daskin [14]

The vehicle fleet consists of 2 homogeneous electric vehicles, each with a passenger capacity of 3 ($Q^k = 3$) and a battery capacity of 16 kWh ($B_{cap}^k = 16 \, kWh$). The vehicles' energy consumption is 6.9 kWh/hr, so each vehicle can travel 110 km with a fully charged battery. The maximum allowable time that vehicle k can be used is considered to be 8 hours, as drivers cannot work more than 8 hours a day ($T^k = 8hr$). The minimum and maximum allowable proportions of the battery capacity are 0.3 and 0.9, respectively ($B_{Min} = 0.3, B_{Max} = 0.9$). Moreover, vehicles have to charge at least 0.7 of their battery capacity each time they charge their batteries ($C_{Min} = 0.7$).

There are two charging stations available at the initial depot and throughout the network. The one at the initial depot is a moderate station with no queue, and the other one is a slow station with a time-dependent queue. The amount of time each vehicle spends at a charging station depends on the time window of the subsequent node it has to visit. The charging rate of vehicles at each station is determined based on the values shown in Table 1.

The objective function weight factors ($w_1, w_2, w_3$) are set to {0.75, 0.25, 0.75}. That is, the total vehicle travel time and penalty of late arrivals at the final depot account for nearly 80% of the objective function score, and the deviation from passengers' desired time windows accounts for 20%. This is because the number of demands is relatively low in our problem, so more emphasis has been placed on reducing total travel time and drivers' excess ride time.

The pick-up nodes (red dots) and drop-off nodes (green dots), as well as the public queued charging station (blue dot), are depicted in Figure 5. The initial and final depots, along with the private queueless charging station, are located at the center of the area and are also indicated with blue dots. The coordinates of these nodes, along with the variation in the number of passengers at each node ($q_i$), the desired time windows, and the service time at each node ($dstop_i$), are indicated in Table 5.

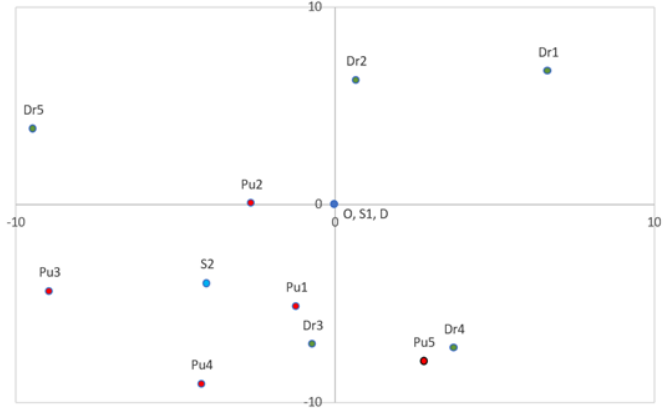

**Figure 5-** The $20 \times 20 \ km^2$ network

**Table5**

Description of different nodes through out the network:

| Symbol | Definition | X-coordinates | Y-coordinates | $q_i$ | Desired time window | $dstop_i$ (min) |
|---|---|---|---|---|---|---|
| $pu_1$ | Passenger 1's prick-up node | -1.198 | -5.164 | 1 | [8:00, 15:00] | 3 |
| $pu_2$ | Passenger 2's prick-up node | -2.61 | 0.039 | 1 | [8:32, 8:47] | 3 |
| $pu_3$ | Passenger 3's prick-up node | -8.938 | -4.388 | 1 | [8:14, 8:29] | 3 |
| $pu_4$ | Passenger 4's prick-up node | -4.172 | -9.096 | 1 | [11:18, 11:33] | 3 |
| $pu_5$ | Passenger 5's prick-up node | 2.792 | -7.944 | 1 | [11:00, 11:15] | 3 |
| $dr_1$ | Passenger 1's drop-off node | 6.687 | 6.731 | -1 | [14:42, 14:57] | 3 |
| $dr_2$ | Passenger 2's drop-off node | 0.673 | 6.283 | -1 | [8:00, 15:00] | 3 |
| $dr_3$ | Passenger 3's drop-off node | -0.694 | -7.098 | -1 | [8:00, 15:00] | 3 |
| $dr_4$ | Passenger 4's drop-off node | 3.736 | -7.269 | -1 | [8:00, 15:00] | 3 |
| $dr_5$ | Passenger 5s drop-off node | -9.45 | 3.792 | -1 | [8:00, 15:00] | 3 |
| $s_1$ | Private charging station | 0 | 0 | 0 | [8:00, 16:00] | 0 |
| $s_2$ | Public charging station | -4 | -4 | 0 | [8:00, 16:00] | 0 |
| $o(0)$ | Initial depot | 0 | 0 | 0 | [8:00, 16:00] | 0 |
| $d(11)$ | Final depot | 0 | 0 | 0 | [8:00, 16:00][1] | 0 |

---

[1] The time windows mentioned for $s_{11}, s_{21}, o_1, o_2$, and $d$ nodes are, in fact, their operating hours.

Furthermore, the Euclidean distances between nodes are calculated and displayed in Table 6. These values are divided by 6.9 $\frac{Km}{KWh}$ and the speed value of the step function ($\frac{V_{ij}^m}{60}$) to determine the energy consumption parameter ($e_{ij}$) and travel time ($t_{i,j}^m$) for each link.

**Table6**

Euclidean distances (km) between nodes:

| i \ j | $o$ | $pu_1$ | $pu_2$ | $pu_3$ | $pu_4$ | $pu_5$ | $dr_1$ | $dr_2$ | $dr_3$ | $dr_4$ | $dr_5$ | $d$ | $s_{11}$ | $s_{21}$ |
|---|---|---|---|---|---|---|---|---|---|---|---|---|---|---|
| $o$ | 0 | 5.3 | 2.6 | 10.0 | 10.0 | 8.4 | 9.5 | 6.3 | 7.1 | 8.2 | 10.2 | 0 | 0 | 5.7 |
| $pu_1$ | 5.3 | 0 | 5.4 | 10.2 | 4.9 | 4.9 | 14.3 | 11.6 | 2.0 | 5.4 | 12.2 | 5.3 | 5.3 | 3.0 |
| $pu_2$ | 2.6 | 5.4 | 0 | 12.4 | 9.3 | 9.6 | 11.5 | 7.1 | 7.4 | 9.7 | 7.80 | 2.6 | 2.6 | 4.3 |
| $pu_3$ | 10.0 | 10.2 | 12.4 | 0 | 14.0 | 7.1 | 11.3 | 13.5 | 10.0 | 5.9 | 20.1 | 10.0 | 10.0 | 12.9 |
| $pu_4$ | 10.0 | 4.9 | 9.27 | 13.4 | 0 | 7.1 | 19.2 | 16.1 | 4.0 | 8.1 | 13.4 | 10.0 | 10.0 | 5.1 |
| $pu_5$ | 8.4 | 4.9 | 9.6 | 7.1 | 7.1 | 0 | 15.2 | 14.4 | 3.6 | 1.2 | 17.0 | 8.4 | 8.4 | 7.9 |
| $dr_1$ | 9.5 | 14.3 | 11.5 | 11.3 | 19.2 | 15.2 | 0 | 6.0 | 15.6 | 14.3 | 16.4 | 9.5 | 9.5 | 15.1 |
| $dr_2$ | 6.3 | 11.6 | 7.1 | 13.5 | 16.1 | 14.4 | 6.0 | 0 | 13.5 | 13.9 | 10.4 | 6.3 | 6.3 | 11.3 |
| $dr_3$ | 7.1 | 2 | 7.4 | 10.0 | 4.0 | 3.6 | 15.7 | 13.05 | 0 | 4.5 | 14.0 | 7.1 | 7.1 | 4.5 |
| $dr_4$ | 8.2 | 5.4 | 9.7 | 5.9 | 8.1 | 1.2 | 14.3 | 13.9 | 4.5 | 0 | 17.2 | 8.2 | 8.2 | 8.4 |
| $dr_5$ | 10.2 | 12.2 | 7.8 | 20.1 | 13.9 | 17.0 | 16.4 | 10.4 | 14.0 | 17.2 | 0 | 10.2 | 10.2 | 9.5 |
| $d$ | 0 | 5.3 | 2.6 | 10.0 | 10.0 | 8.4 | 9.5 | 6.3 | 7.1 | 8.2 | 10.2 | 0 | 0 | 5.7 |
| $s_{11}$ | 0 | 5.3 | 2.6 | 10.0 | 10.0 | 8.4 | 9.5 | 6.3 | 7.1 | 8.2 | 10.2 | 0 | 0 | 5.7 |
| $s_{21}$ | 5.7 | 3.0 | 4.3 | 13.0 | 5.1 | 7.9 | 15.1 | 11.3 | 4.5 | 8.4 | 9.5 | 5.7 | 5.7 | 0 |

We inserted all the parameters into an Excel sheet, formulated the model in the GAMS software, and then solved it using the CPLEX software on a personal computer with a 16 GB RAM and 2.30 GHz Core i7 CPU to test all of the model's features. To fully capture the vehicles' performance under different initial battery levels, traffic congestion, and passenger utility, we defined three scenarios.

### 4.1 The results of the first scenario

In the first scenario, the basic model introduced in Section 3, with 1820 variables, was solved in 830 seconds (nearly 14 minutes). The results indicate that the two vehicles could satisfy all the demand. The optimal planned routes of these vehicles are depicted in Figures 6 and 7.

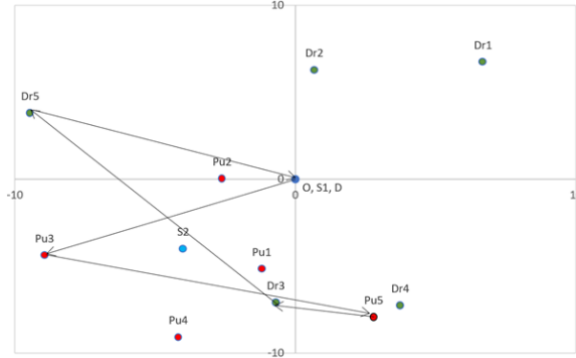

1
2  **Figure 6-** First vehicle's route in the first scenario
3

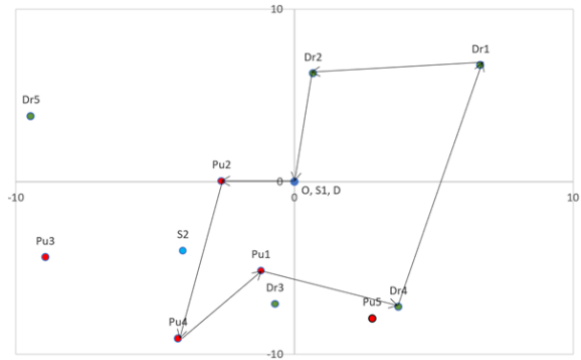

4
5  **Figure 7-** Second vehicle's route in the first scenario
6
7        According to the results of the model mentioned in Table 7, the first vehicle leaves the initial
8    depot at 8:00 ($o$) to pick up the third passenger at 8:16 ($pu_3$). Then, it heads to the fifth passenger's
9    pick-up node and reaches there by 11:00 ($pu_5$). Then the third and fifth passengers were dropped off
10   at their destinations at 13:17 ($dr_5$) and 13:41 ($dr_5$), respectively. Eventually, the first vehicle went to
11   the final depot at 13:59 after traversing 44.8km without charging its battery and spending nearly 7
12   hours in the network, providing service without any deviation from the announced time windows.
13       **Table 7**
14   Detailed results of the model for the first scenario:

| Vehicle | Served node | Service time | Desired time window |
|---------|-------------|--------------|---------------------|
| 1 | $o$ | 8:00 | [8:00, 16:00] |
| 1 | $pu_3$ | 8:16 | [8:14, 8:29] |
| 1 | $pu_5$ | 11:00 | [11:00, 11:15] |
| 1 | $dr_3$ | 13:17 | [8:00, 15:00] |
| 1 | $dr_5$ | 13:41 | [8:00, 15:00] |
| 1 | $d$ | 13:59 | [8:00, 16:00] |
| 2 | $o$ | 8:27 | [8:00, 16:00] |
| 2 | $pu_2$ | 8:32 | [8:14, 8:29] |
| 2 | $pu_4$ | 11:18 | [11:00, 11:15] |
| 2 | $pu_1$ | 13:17 | [11:00, 11:15] |

| 2 | $dr_4$ | 13:28 | [8:00, 15:00] |
| 2 | $dr_1$ | 14:42 | [8:00, 15:00] |
| 2 | $dr_2$ | 14:54 | [8:00, 15:00] |
| 2 | $d$ | 15:07 | [8:00, 16:00] |

Vehicle 2 picked up the second, fourth, and first passengers at 8:32, 11:18, and 13:17, respectively. In this scenario, Vehicle 2 could have left the initial depot at 8:00. However, it chose to depart at 8:27 instead, to avoid any deviation from the second passenger's desired time window. This vehicle dropped off the fourth, first, and second passengers at 13:28, 14:32, and 14:54, respectively, before reaching the final depot at 15:07. It covered a distance of 48.85 km without charging its battery.

The impact of traffic congestion on the results is readily observable. For instance, the first vehicle was able to traverse a distance of 10 km from the initial depot to the third passenger's pick-up node in just 16 minutes during the morning time period ([8:00, 9:20]). However, during peak hours, the same vehicle took approximately 2 hours to cover a shorter distance of 3.6 km between the locations $pu_5$ and $dr_3$. Furthermore, the third and second passengers were served first as these two had the earliest desired pick-up time window. The amount of in-vehicle time that passengers 1 to 5 spent in this scenario is 1.4, 6.4, 5, 2.2, and 2.7 hours, respectively.

To further demonstrate the effectiveness of shared mobility, we computed the total travel time and distance in the network for the private mode and compared the results. Although the total travel time in shared mobility has increased by about 6%, from 84.63 $\frac{min}{veh}$ in the private mode to 89.32 $\frac{min}{veh}$ in the shared mode, the total distance covered shows a 20% drop, from 56.42 $\frac{km}{veh}$ in the private mode to 46.83 $\frac{km}{veh}$ in the shared mode. This indicates the effectiveness of shared mobility in reducing traffic congestion, which in turn could provide economic and environmental benefits. However, the amount of in-vehicle time for some passengers, e.g., passenger 2, is very high, which could discourage passengers in the real world from adopting shared mobility. To address this issue and increase passenger utility, we introduced a second scenario where the in-vehicle time of passengers has been confined by adding a constraint to the basic formulation.

### 4.2 The results of the second scenario

In this scenario, Equation (53) was added to the basic formulation to limit the passengers' in-vehicle time to twice the private travel time plus 10 minutes.

$$U_{n+i}^k - U_i^k \leq 2 * \sum_j \sum_m t_{i,n+i}^m x_{i,j}^{k,m} + 10 \text{ (min)}; \ \forall \ i \in P, k \in K \tag{53}$$

The computational time for this scenario was 130 seconds. The vehicles' routes are depicted in Figures 8 and 9. These results show a significant deviation from the results of the first scenario. For instance, given the pick-up locations and desired time windows, the first vehicle could have picked up the fifth passenger after picking up passenger 2 at 8:32. However, this would have violated Equation (53). Therefore, it dropped off passenger 2 at 8:45 before serving passenger 5. Further details are provided in Table 8.

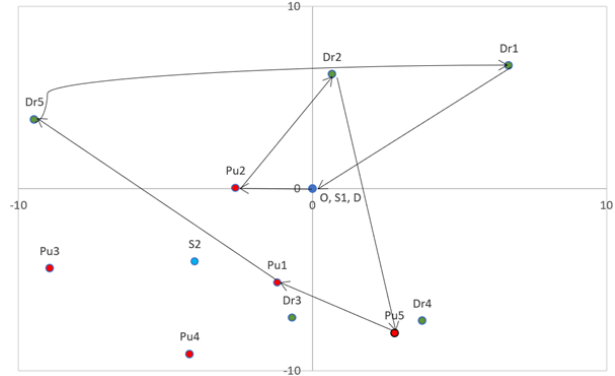

**Figure 8-** First vehicle's route in the second scenario

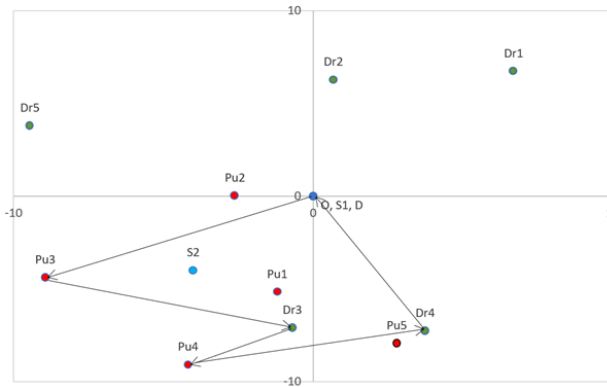

**Figure 9-** Second vehicle's route in the second scenario

**Table8**

Detailed results of the model for the second scenario:

| Vehicle | Served node | Service time | Desired time window |
|---------|-------------|--------------|---------------------|
| 1 | $o$ | 8:27 | [8:00, 16:00] |
| 1 | $pu_2$ | 8:32 | [8:32, 8:47] |
| 1 | $dr_2$ | 8:45 | [8:00, 15:00] |
| 1 | $pu_5$ | 11:00 | [11:00, 11:15] |
| 1 | $pu_1$ | 13:49 | [8:00, 15:00] |
| 1 | $dr_5$ | 14:10 | [8:00, 15:00] |
| 1 | $dr_1$ | 14:42 | [14:42, 14:57] |
| 1 | $d$ | 14:59 | [8:00, 16:00] |
| 2 | $o$ | 8:00 | [8:00, 16:00] |
| 2 | $pu_3$ | 8:16 | [8:14, 8:29] |
| 2 | $dr_3$ | 8:34 | [8:00, 15:00] |
| 2 | $pu_4$ | 11:29 | [11:18, 11:33] |
| 2 | $dr_4$ | 13:17 | [8:00, 15:00] |
| 2 | $d$ | 13:32 | [8:00, 16:00] |

Adopting the second scenario has reduced the in-vehicle time for passengers 1 to 4, albeit at the cost of an increase in both the total distance covered and total travel time. To elaborate, the in-vehicle time for passengers 1 to 5 is 0.9, 0.2, 0.3, 1.8, and 3.2 hours, respectively. It's worth noting that the fifth passenger's desired time window is [11:00, 11:15], which falls within the day's peak period. Due to the time-dependency of the model, this passenger's in-vehicle time did not decrease compared to the first scenario, unlike the other passengers who experienced significant reductions in in-vehicle time. It should be highlighted that passenger 5 would have experienced a considerable in-vehicle time even in private mode, as they are traveling during the day's peak period. The network's performance under the second scenario, in terms of total travel time ($\frac{min}{veh}$), is 109.61 $\frac{min}{veh}$ in the shared mode and 84.63 $\frac{min}{veh}$ in the private mode. The total distance covered (km/veh) is 53.64 $\frac{km}{veh}$ in the shared mode and 56.42 $\frac{km}{veh}$ in the private mode.

### 4.3 The results of the third scenario

In previous scenarios, vehicles did not have to charge their batteries during their trips as they left the initial depot with the maximum battery level (16 kWh). Hence, we defined the third scenario where they commence their trips with 50% battery capacity (8 kWh) to analyze the changes in the system's performance. The routes of the first and second vehicles are depicted in Figures 10 and 11.

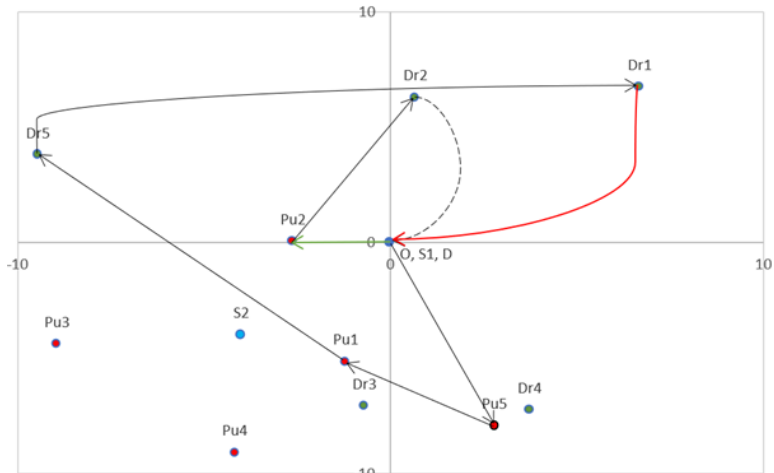

**Figure** First vehicle's route in the third scenario

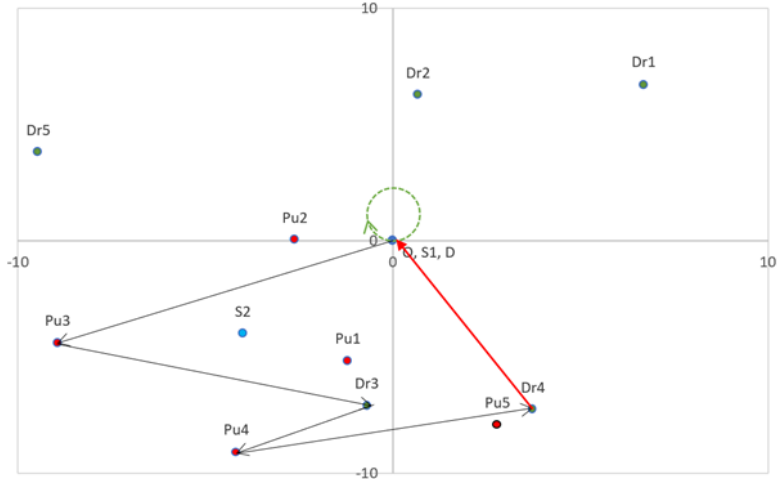

**Figure 11-** Second vehicle's route in the third scenario

The service plan for the first vehicle is similar to its plan in the second scenario, with only one difference. Its battery level reached 6.56 kWh after dropping off the second passenger, and if it decided to pick up passenger 5, its battery level would have dropped to 4.5 kWh, which is below the minimum allowable battery level ($0.3 \times 16 = 4.8$ kWh). Therefore, it heads to the private charging station at the initial depot with a battery level of 5.684 kWh, and leaves after 12.517 minutes with a battery level of 12.242 kWh. Its nonlinear charging function, along with the linear approximation functions, are depicted in Figure 12.

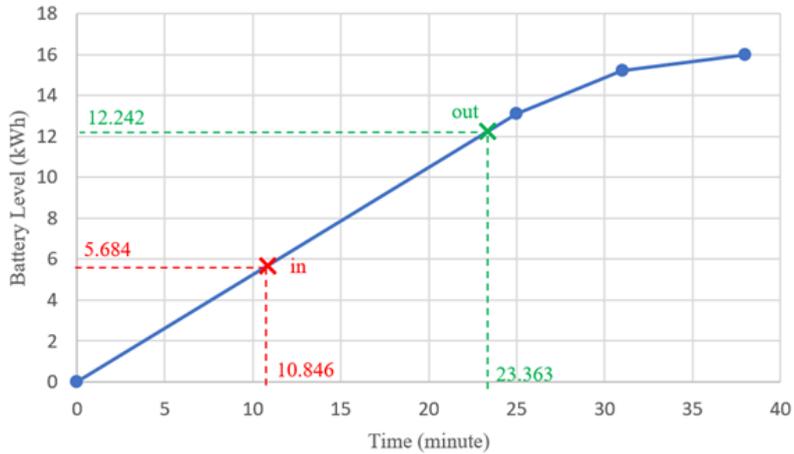

**Figure 12-** First vehicle's nonlinear charging function in the third scenario

The service sequence of the second vehicle is also similar to its plan in the second scenario. However, it has decided to stay at the initial depot for 6.107 minutes and charge its battery to 11.2 kWh before providing service to the third and fourth passengers. Its charging function is depicted in Figure 13. Further details of the routing plan is provided in Table 9.

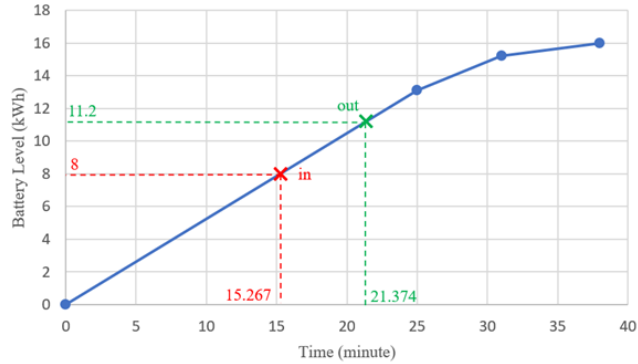

**Figure 13-** Second vehicle's nonlinear charging function in the third scenario

**Table9**
Detailed results of the model for the third scenario:

| Vehicle | Served node | Service time | Desired time window |
|---------|-------------|--------------|---------------------|
| 1 | $o$ | 8:27 | [8:00, 16:00] |
| 1 | $pu_2$ | 8:32 | [8:32, 8:47] |
| 1 | $dr_2$ | 8:45 | [8:00, 15:00] |
| 1 | $s_1$ | 9:20 | [8:00, 16:00] |
| 1 | $pu_5$ | 11:00 | [11:00, 11:15] |
| 1 | $pu_1$ | 13:49 | [8:00, 15:00] |
| 1 | $dr_5$ | 14:10 | [8:00, 15:00] |
| 1 | $dr_1$ | 14:42 | [14:42, 14:57] |
| 1 | $d$ | 14:59 | [8:00, 16:00] |
| 2 | $o$ | 8:00 | [8:00, 16:00] |
| 2 | $s_1$ | 8:01 | [8:00, 8:16] |
| 2 | $pu_3$ | 8:24 | [8:14, 8:29] |
| 2 | $dr_3$ | 8:42 | [8:00, 15:00] |
| 2 | $pu_4$ | 11:29 | [11:18, 11:33] |
| 2 | $dr_4$ | 13:17 | [8:00, 15:00] |
| 2 | $d$ | 13:32 | [8:00, 16:00] |

## 5. Conclusion

In this study, we proposed an integer programming model to optimize the routing plan of an electric vehicle fleet in a time-dependent dial-a-ride transit system. We attempted to incorporate traffic congestion into the network by adopting a time-dependent vehicle speed function, a nonlinear charging function, considerations for the vehicle's battery health preservation, and different charging stations with time-dependent waiting times to make the results more applicable in the real world.

The results of the case study enabled us to gain several insights such as: (i) the battery level of electric vehicles could be managed by enabling them to charge only up to a minimum value, without being forced to fully charge their batteries; (ii) considering an upper limit for the passenger's in-vehicle time could increase the passenger's utility by preventing high travel times without making a significant impact on the fleets' total distance traveled and total travel time, hence, encouraging more people to use shared mobility; (iii) a well-planned shared transit system could effectively enhance the network's performance by reducing the total distance traveled without significantly increasing the total travel time. This, in turn, could help alleviate traffic congestion and its problematic consequences, such as air pollution and operational costs in urban areas.

Although the results of the case study were plausible, several suggestions can enhance the model's performance and practicality such as: (i) considering a continuous speed-time function instead of a step function could capture the time dependency of the network more precisely. This could be accomplished by considering the demand assignment in the formulation; (ii) solving the model for large-scale networks with high demand using CPLEX is burdensome as the computational time of the basic model for the small-scale network of the case study with only 5 passenger demand took 14 minutes. Therefore, future works could focus on developing heuristic solutions to this problem; (iii) the results of the second and third scenarios in the case study indicated that the routing plan produced by the model could result in some private trips. To confine the number of private trips, another term could be added to the objective function to minimize the private trips as much as possible; (iv) adjusting the model to consider dynamic demand where passengers could cancel their trips or add another trip demand in the middle of the day; (v) considering a pricing policy to control traffic congestion during the peak period and prevent unnecessary trips or trip demands with very narrow desired time windows.

**Ethical Statement**

The authors declare that they have no known competing financial interests or personal relationships that could have appeared to influence the work reported in this paper.

**Author Contribution Statement**

The authors confirm contribution to the paper as follows: study conception and design: A. Yazdiani, S. Bafandkar, A. Elmi, Y. Shafahi,; analysis and interpretation of results: A. Yazdiani, Y. Shafahi; draft manuscript preparation: S. Bafandkar. All authors reviewed the results and approved the final version of the manuscript.